\documentclass[11pt]{article}
\usepackage{amssymb,amsfonts,amsmath,curves}
\usepackage{epsfig}
\usepackage{graphicx}
\usepackage{color}
\usepackage{anysize}

\makeatletter\@addtoreset{equation}{section}
\makeatother

\newtheorem{defin}{Definition}[section]

\newtheorem{teo}{Theorem}[section]

\newtheorem{con}{Conjecture}[section]
\newtheorem{cond}{Condition}
\newtheorem{conj}{Conjecture}[section]
\newtheorem{prop}[teo]{Proposition}
\newtheorem{lem}{Lemma}[section]
\newtheorem{rmk}[teo]{Remark}
\newtheorem{cor}{Corollary}[section]

\newcommand{\be}{\begin{equation}}
\newcommand{\ee}{\end{equation}}
\newcommand{\ben}{\begin{eqnarray}}
\newcommand{\benn}{\begin{eqnarray*}}
\newcommand{\een}{\end{eqnarray}}
\newcommand{\eenn}{\end{eqnarray*}}
\newcommand{\bp}{\begin{prop}}
\newcommand{\ep}{\end{prop}}
\newcommand{\bt}{\begin{teo}}
\newcommand{\et}{\end{teo}}
\newcommand{\bcor}{\begin{cor}}
\newcommand{\ecor}{\end{cor}}
\newcommand{\bcon}{\begin{con}}
\newcommand{\econ}{\end{con}}
\newcommand{\bcond}{\begin{cond}}
\newcommand{\econd}{\end{cond}}
\newcommand{\bconj}{\begin{conj}}
\newcommand{\br}{\begin{rmk}}
\newcommand{\er}{\end{rmk}}
\newcommand{\bl}{\begin{lem}}
\newcommand{\el}{\end{lem}}
\newcommand{\bit}{\begin{itemize}}
\newcommand{\eit}{\end{itemize}}
\newcommand{\bd}{\begin{defin}}
\newcommand{\ed}{\end{defin}}
\newcommand{\bpr}{\begin{proof}}
\newcommand{\epr}{\end{proof}}

\newenvironment{proof}{\noindent {\em Proof}.\,\,}{\hspace*{\fill}$\halmos$\medskip}

\newcommand{\halmos}{\rule{1ex}{1.4ex}}
\def \qed {{\hspace*{\fill}$\halmos$\medskip}}

\newcommand{\Z}{{\mathbb Z}}
\newcommand{\R}{{\mathbb R}}
\newcommand{\E}{{\mathbb E}}

\renewcommand{\P}{{\mathbb P}}
\newcommand{\N}{{\mathbb N}}

\marginsize{1.8cm}{1.8cm}{1.0cm}{2.2cm}

\begin{document}

\title{A quenched limit theorem for the local time of random walks on $\Z^2$}
\author{J\"urgen G\"artner$^{\, 1,2}$ \and Rongfeng Sun$^{\, 1,3}$}
\maketitle

\footnotetext[1]{MA 7-5, Fakult\"at II -- Institut f\"ur Mathematik, TU Berlin, Stra\ss e des 17.\ Juni 136,
10623 Berlin.}
\footnotetext[2]{\it jg@math.tu-berlin.de}
\footnotetext[3]{\it sun@math.tu-berlin.de}

\vspace{-16pt}

\begin{abstract}
Let $X$ and $Y$ be two independent random walks on $\Z^2$ with zero mean and finite variances, and let $L_t(X,Y)$ 
be the local time of $X-Y$ at the origin at time $t$. We show that almost surely with respect to $Y$, $L_t(X,Y)/\log t$
conditioned on $Y$ converges in distribution to an exponential random variable with the same mean as the
distributional limit of $L_t(X,Y)/\log t$ without conditioning. This question arises naturally from the study of 
the parabolic Anderson model with a single moving catalyst, which is closely related to a pinning model.

\medskip


\end{abstract}

\section{Introduction}

It is a classical result dating back to Erd\"os and Taylor \cite{ET60} that, for a simple random walk on $\Z^2$,
if $L_t$ denotes its local time at the origin at time $t$, then $L_t/\log t$ converges in distribution to
an exponential random variable as $t\to\infty$. With a change of parameter for the exponential random variable, the 
same result holds for general zero mean finite variance random walks on $\Z^2$. More precisely, if $X$ is either
a discrete or a continuous time random walk on $\Z^2$ with zero mean, finite variance, and one-step increment
distribution $p(\cdot)$, then its covariance matrix is defined by 
\be\label{Q}
Q_{ij}= \sum_{x\in\Z^2} p(x) x_ix_j, \qquad  1\leq i, j\leq 2. 
\ee
Let $L_t(X)=\sum_{i=0}^t \delta_0(X_i)$ if $X$ is a discrete time random walk, and let $L_t(X)=\int_0^t \delta_0(X_s)ds$
if $X$ is a continuous time random walk. Then the classical Erd\"os--Taylor result states that

\bt\label{T:explaw}{\bf[Erd\"os--Taylor]}
Let $X$ be an irreducible zero mean finite variance random walk on $\Z^2$ with covariance matrix $Q$ starting at the
origin. Let $r$ denote the jump rate of $X$ if it is a continuous time random walk, and set $r=1$ otherwise.
Then as $t\to\infty$, $\E\left[\Big(\frac{2\pi r\sqrt{\det Q}L_t}{\log t}\Big)^k\right] \to  k!$ for each $k\in\N$, and 
$\frac{2\pi r \sqrt{\det Q} L_t}{\log t}$ converges in distribution to a mean 1 exponential random variable.
\et
{\bf Remark.} If $X$ is not irreducible, but is still truly two-dimensional, then the sublattice in 
$\Z^2$ which $X$ visits with positive probabilty can be mapped {\it linearly and 
bijectively} to 
$\Z^2$ (see P1 in Sec.\ 7 and P5 in Sec.\ 2 of Spitzer \cite{S76}). Theorem \ref{T:explaw} can then be applied 
to the image random walk. 
\medskip

If $X$ and $Y$ are two independent, but not necessarily identically distributed, irreducible zero mean finite
variance random walks on $\Z^2$ such that $X-Y$ is also irreducible, then Theorem \ref{T:explaw} applies to $X-Y$. 
This can be regarded as an averaged limit theorem for the local time $L_t(X,Y):= L_t(X-Y)$, where $Y$ plays the role 
of the random environment. The objective of this paper is to obtain a quenched limit 
theorem for $L_t(X,Y)$, i.e., a limit theorem 
for $L_t(X,Y)$ conditioned on $Y$.

For future reference, let $\P^X_x(\cdot)$ denote probability w.r.t.\ the random walk $X$ starting from $x$,
and let $\E^X_x[\cdot]$ denote the corresponding expectation.

\bt\label{T:2rw}{\bf[Quenched exponential law]}
Let $X$ and $Y$ be independent irreducible zero mean finite variance random walks on $\Z^2$ starting from the origin,
such that $Z:=X-Y$ is also irreducible. Let $Q$ be the covariance matrix of $Z$. Let $\kappa>0$ and $\rho>0$ denote 
the respective jump rates of $X$ and $Y$ if they are continuous time random walks, and set $\kappa+\rho=1$ if they 
are discrete time random walks. Then almost surely with respect to $Y$, as $t\to\infty$, 
$\E^X_0\big[(\frac{2\pi (\kappa+\rho)\sqrt{\det Q}L_t(X,Y)}{\log t})^k |Y\big]\to k!$ for each $k\in\N$, and 
$\frac{2\pi(\kappa+\rho)\sqrt{\det Q} L_t(X,Y)}{\log t}$ conditioned on $Y$ converges in distribution to a mean 
1 exponential random variable.
\et
{\bf Remark.} If $Z$ is reducible, e.g., when $X$ and $Y$ are discrete time 
simple random walks on $\Z^2$, then $Q$ needs to be replaced by the covariance matrix of an
image random walk, namely, the random walk obtained from $Z$ after one applies the {\it linear} map
which maps the set of sites in $\Z^2$ that $Z$ visits with positive probability {\it bijectively} to $\Z^2$.
\bigskip

\noindent
{\bf Remark.} The analogue of Theorem \ref{T:2rw} fails for dimensions $d\neq 2$. Consider the discrete time case.
For $d\geq 3$, by the transience of the random walk $X-Y$, a.s.\ w.r.t.\ $X$ and $Y$, $L_n(X,Y)$ increases to a random 
constant $L_\infty(X,Y)$ as $n\to\infty$. With respect to the joint law of $X$ and $Y$, $L_\infty(X,Y)$ is 
geometrically distributed; however conditioned on $Y$, the law of $L_\infty(X,Y)$ clearly depends sensitively on the 
realization of $Y$. For $d=1$, the correct scaling for $L_n(X,Y)$ is $n^{-1/2}$. Under diffusive scaling, $(X,Y)$ 
converges in law to a pair of independent Brownian motions $(B_1, B_2)$, while up to a constant factor, 
$L_n(X,Y)/\sqrt{n}$ converges in law to the collision local time $\bar L_1(B_1, B_2)$ between $B_1$ and $B_2$ up to 
time $1$. Thus as random probability distributions, the law of $L_n(X,Y)/\sqrt{n}$ conditioned on $Y$ is expected
to converge to the law of $\bar L_1(B_1, B_2)$ conditioned on $B_2$. However, such a convergence will only take place 
in probability instead of a.s., because the law of $L_n(X,Y)/\sqrt{n}$ conditioned on $Y$ depends
sensitively on the rescaled path $(Y_i/\sqrt{n})_{0\leq i\leq n}$, which a.s.\ does not converge as $n\to\infty$.
We will not pursue the $d=1$ case in this paper.
\bigskip

Our original motivation for the study of the law of $L_t(X,Y)$ conditioned on $Y$ stems from the parabolic Anderson model
where the random medium consists of a single moving catalyst:
\be \label{PAMs}
\begin{aligned}
\frac{\partial}{\partial t} u(t, x) &= \kappa \Delta u(t, x) + \gamma\delta_{Y_t}(x) \, u(t, x),  \\
u(0,x) &= 1,
\end{aligned}
\qquad \qquad x\in \Z^d,\ t\geq 0,
\ee
where $\kappa\geq 0$, $\gamma\in\R$, $\Delta f(x) = \frac{1}{2d} \sum_{\Vert y-x\Vert =1} (f(y)-f(x))$ is the
discrete Laplacian on $\Z^d$, and $Y_t$ is a simple random walk on $\Z^d$ with jump rate $\rho\geq 0$. By the
Feynman-Kac representation,
\be\label{FK}
u(t, x) = \E^X_x\left[\exp \left\{\gamma \int_0^t \delta_0(X_s-Y_{t-s})ds\right\}\right],
\ee
where $X$ is a simple random walk on $\Z^d$ with jump rate $\kappa$ and starting from $x$. Note that if not for the 
time reversal of $Y$ in (\ref{FK}), the exponent in (\ref{FK}) would be exactly 
$\gamma L_t(X,Y)$. 

The annealed Lyapunov exponents $\lambda_k = \lim_{t\to\infty} t^{-1} \log \E^Y_0[u(t,0)^k]$, $k\in\N$, were studied 
by G\"artner and Heydenreich in \cite{GH06}. For the quenched Lyapunov 
exponent $\lambda= \lim_{t\to\infty} t^{-1} \log u(t,0)$, we can replace 
$u(t,0)$ by $\underline u_{0,t}$, where
\be
\underline{u}_{s,t} = 
\E^X_{Y_t}\left[\exp \left\{\gamma \int_0^{t-s} \delta_0(X_a-Y_{t-a})da\right\} 1_{\{X_{t-s}=Y_s\}}\right], \qquad \qquad 0\leq
s <t.
\ee
It turns out that $\lambda = \lim_{t\to\infty} t^{-1}\log u(t,0)=\lim_{t\to\infty} t^{-1} \log 
\underline{u}_{0,t}$. By the superadditive ergodic
theorem applied to
$\log\underline{u}_{s,t}$, it can be shown that 
\be\label{lam}
\lambda = \sup_{t>0} \frac{1}{t}\, \E^Y_0[\log\underline{u}_{0,t}] 
= \sup_{t>0} \frac{1}{t}\, \E^Y_0\Big[\log \E^X_0\big[e^{\gamma L_t(X,Y)}1_{\{X_t=Y_t\}}\big]\Big],
\ee
where we have reversed time for $Y$ in the second equality.
There exists a critical $\gamma_c\in\R$ such that $\lambda=0$ if $\gamma\leq \gamma_c$, and $\lambda>0$ if
$\gamma >\gamma_c$. It can be shown that $\gamma_c=0$ in dimensions $d=1,2$, and $\gamma_c>0$ in $d\geq 3$.
The proof of $\gamma_c=0$ in $d=2$ is the most subtle one, and the only proof we know of at the moment 
uses the representation (\ref{lam}) and Theorem \ref{T:2rw}. The details are contained in Birkner and Sun 
\cite{BS08}. 

A closely related model where the conditional law of $L_t(X,Y)$ arises naturally is a pinning model. More precisely,
we define a change of measure from the random walk path measure $P$ on $(X_s)_{0\leq s\leq t}$ with 
Radon-Nikodyn derivative 
\be
\frac{d P^\gamma_{t, Y}}{d P} = \frac{e^{\gamma L_t(X,Y)}}{Z^\gamma_{t,Y}},
\ee 
where $Z^\gamma_{t,Y} = \E^X_0[e^{\gamma L_t(X,Y)}]$ is the normalizing constant. With respect to the measure
$P^\gamma_{t,Y}$, the random walk $X$ prefers to be at the same location as $Y$ when $\gamma>0$. This model
exhibits a localization-delocalization transition. Namely, there exists a critical $\gamma_c\in\R$
such that if $\gamma< \gamma_c$, then for typical $Y$ and typical $X$ w.r.t.\ $P^\gamma_{t,Y}$, $X$ and $Y$
spend negligible fraction of time together; while if $\gamma >\gamma_c$, then for typical $Y$ and typical
$X$ w.r.t.\ $P^\gamma_{t,Y}$, $X$ and $Y$ spend positive fraction of time together. By the same argument as for the 
parabolic Anderson model (\ref{PAMs}), it can be shown that $\lim_{t\to\infty} t^{-1}\log Z^\gamma_{t,Y}$, the 
so-called free energy, exists almost surely and equals $\lambda$ in (\ref{lam}) (see \cite{BS08} for details).
This implies that $\gamma_c=0$ in $d=1, 2$, and $\gamma_c>0$ in $d\geq 3$. For more on pinning models in general, 
see Giacomin \cite{G07}.

Another model where the conditional law of $L_t(X,Y)$ appears is the directed polymer model in random environment.
See Birkner \cite{B04} for a sufficient condition for weak disorder which is formulated in terms of the law of 
$L_t(X,Y)$ conditioned on $Y$.

The exponential law arises in many different contexts in the study of the local time of two-dimensional random walks. 
Another interesting instance is a result by \v{C}ern\'y \cite{C07} that, almost surely with respect to the path of a 
non-degenerate zero mean finite variance random walk on $\Z^2$, as 
$t\to\infty$, the law of the local time at time $t$ sampled 
uniformly among all sites visited by the walk up to time $t$, and rescaled by a factor of $1/\log t$, converges to the 
law of an exponential random variable.
\bigskip

To end the introduction, we propose an interesting open problem.
\medskip

\noindent{\bf Open Problem:}
Fix $k\geq 1$. Let $X$, $Y_1$, \ldots,$Y_k$ be independent irreducible zero mean finite variance random walks on $\Z^2$ 
starting from the origin, such that $Z_i:=X-Y_i$, $1\leq i\leq k$, are all irreducible. Is it true that as $t\to\infty$,
a.s.\ w.r.t.\ $Y_1, \cdots, Y_k$, $\big(\frac{L_t(X,0)}{\log t}, \frac{L_t(X,Y_1)}{\log t}, \ldots, \frac{L_t(X,Y_k)}{\log
  t}\big)$ conditioned on $Y_1, \cdots, Y_k$ converge in distribution to $k+1$ {\it independent} exponential random variables?
\medskip

Preliminary calculations of expressions of the form $\E^{Y_1}_0\big[\E^X_0\big[\frac{L_t(X,Y_1)}{\log t} e^{-\frac{\gamma
    L_t(X,0)}{\log t}}\big]\big]$, assuming the quantity inside $\E^{Y_1}_0[\cdot]$ asymptotically self-averages,
favor the affirmative. However, we will not go as far as to formulate it as a conjecture here.

\section{Preliminary Lemmas}

In this section, we prove two lemmas, \ref{L:gradpot} and \ref{L:rwconv}, which we will need to prove
Theorem \ref{T:2rw} in Sec.\ \ref{S3}.
\bl\label{L:gradpot}
Let $Z$ be an irreducible zero mean finite variance random walk on $\Z^2$ with covariance matrix $Q$. Let
$p^Z_n(\cdot)$, resp.\ $p^Z_t(\cdot)$, denote the translation invariant transition probability kernel for the case 
$Z$ is a discrete, resp.\ continuous, time random walk. Then there exists $0<C<\infty$ such that for any 
$x, z_0\in \Z^2$ $($with $p^Z_n(x)p^Z_n(x+z_0)>0$ for some $n\in\N$ in the discrete 
time case$)$, we have
\be\label{A888}
\sum_{n=0}^\infty |p^Z_n(x) - p^Z_n(x+z_0)| \leq C\Vert z_0\Vert\left(\frac{1}{1+\Vert x\Vert}+\frac{1}{1+\Vert x+z_0\Vert}\right),
\ee
where $\Vert\cdot\Vert$ denotes Euclidean norm, and in the continuous time case,
\be\label{A889}
\int_0^\infty |p^Z_t(x) - p^Z_t(x+z_0)| dt \leq C\Vert z_0\Vert\left(\frac{1}{1+\Vert x\Vert} +\frac{1}{1+\Vert x+z_0\Vert}\right).
\ee
\el
{\bf Remark.} The analogue of Lemma \ref{L:gradpot} for random walks on $\Z^d$, $d\geq 3$, is to replace
$1+\Vert x\Vert$ and $1+\Vert x+z_0\Vert$ respectively by $(1+\Vert x\Vert)^{d-1}$ and $(1+\Vert x+z_0\Vert)^{d-1}$
in (\ref{A888}) and (\ref{A889}), which is easily seen if we replace $p^Z_t$ by 
transition densities of Brownian motion.
However, such a result can not hold in general without additional assumptions.
In particular, for $d\geq 4$, we can define a discrete time random walk $p^Z_1(\cdot)$ with 
$p^Z_1(\pm e_i)=1/4d$ for each $1\leq i\leq d$ where
$e_i$ are the unit vectors in $\Z^d$, $p^Z_1(\pm a_ne_1) = Cn^{-2}a_n^{-2}$ for an increasing sequence of
$a_n\in\N$, and $p^Z_1(x)=0$ for all other $x\in\Z^d$. If $a_n$ increases so fast that 
$p^Z_1(\pm a_n e_1) \geq C a_n^{-2-\epsilon}$ for some $\epsilon>0$, then
$|p^Z_1(x+e_1) - p^Z_1(x)|$ already violates the desired decay in $\Vert x\Vert$.

\bl\label{L:rwconv} Let $X$ be an irreducible zero mean finite variance random walk on $\Z^2$. Let $q\in [1,2)$. Then
for all $v\in\Z^2$ and $i\in\N$ $($replace $i\in\N$ by $s\geq 1$ in the continuous time case$)$,
\be\label{A72}
\sum_{x\in\Z^2} \P(X_i=x) \frac{1}{(1+\Vert x-v\Vert)^q} \leq \frac{C_q}{i^{\frac{q}{2}}},
\ee
where $C_q$ is a constant depending only on $q$ and the walk $X$.
\el

To prove Lemma \ref{L:gradpot}, we will use the following expansion form of the local central limit theorem
from Lawler and Limic \cite{LL08} (see Theorem 2.3.8 there for a slightly different formulation).
In \cite{LL08}, this result is stated and proved for discrete time random walks, however it is clear that the 
same proof and result hold for continuous time random walks.

\bt\label{T:LCLT}{\bf [Lawler \& Limic]}
Let $p_n(\cdot)$ be the transition probability kernel of an irreducible aperiodic mean zero random walk on $\Z^d$ 
with finite $(k+1)$-st moment for some integer $k\geq 3$. Let $Q$ be the covariance matrix of the 
random walk. Then 
\be\label{lclt}
p_n(x) = \frac{e^{-\frac{x\cdot Q^{-1}x}{2n}}}{(2\pi n)^{d/2}\sqrt{\det Q}} \left[1+\frac{u_3(x/{\sqrt n})}{\sqrt n}
+\frac{u_4(x/\sqrt{n})}{n}+\cdots +\frac{u_k(x/\sqrt{n})}{n^{(k-2)/2}}\right] + \epsilon_{n,k}(x),
\ee
where there exists $0<c<\infty$ such that
\be\label{uj}
|u_j(z)|\leq c(\Vert z\Vert^j+1),
\ee
and uniformly in $x\in\Z^d$ and $n\in\N$,
\be\label{enk}
|\epsilon_{n,k}(x)| \leq \frac{c}{n^{(d+k-1)/2}}.
\ee
For a rate 1 continuous time random walk, $(\ref{lclt})$--$(\ref{enk})$ hold with $n\in\N$ replaced by $t\geq 1$.
\et
{\bf Proof of Lemma \ref{L:gradpot}.} Initially we only had a proof of Lemma \ref{L:gradpot} for a restricted class
of random walks. Greg Lawler kindly showed us how to extend the result to all irreducible zero mean finite variance
random walks. We present his line of arguments here. Most ingredients can be found in his book with Vlada Limic 
\cite{LL08}. The main idea is to use the {\it finite range coupling} of random walks.

As a remark on notation, since we will not be concerned with the exact values of the constants in our estimates, in
what follows, unless stated otherwise, $c, C, C_1, C_2$, etc, will denote generic constants whose values may change from
line to line.

We only treat the discrete time case. The continuous time case is similar. Without loss of generality,
we may assume that $Z$ is aperiodic, otherwise we can partition $\N$ into periodic subsets and change time scale to
reduce to the aperiodic case. It is not difficult to see that the one-step transition 
kernel $p^Z:=p^Z_1$ allows a decomposition (see Exercise 1.3 in \cite{LL08})
\be\label{ppcouple}
p^Z(x) = \alpha p^{(1)}(x) + (1-\alpha)p^{(2)}(x),
\ee
where $\alpha$ can be chosen in $(0,1/2)$, $p^{(1)}$ is the one-step transitional probability kernel of 
an aperiodic mean zero 
{\it finite range} random walk, and $p^{(2)}$ is the one-step transition probability kernel of an 
aperiodic mean zero finite variance random
walk. Thus a $p^Z$ random walk at each step chooses a jump according to $p^{(1)}$ with probability 
$\alpha$, and
according to $p^{(2)}$ with probability $1-\alpha$. A coupling between two $p^Z$ random walks with different initial
positions is called a {\it finite range coupling} if they choose the same transition kernel
from $\{p^{(1)}, p^{(2)}\}$ at each step, they make the same jumps if $p^{(2)}$ is chosen, and the 
jumps are suitably coupled if $p^{(1)}$ is chosen (see e.g.\ Proposition 2.4.2 and Lemma 2.4.3 in \cite{LL08}).
Let $Q_1$ and $Q_2$ denote respectively the covariance matrices of $p^{(1)}$ and
$p^{(2)}$. If $M_n$ denotes the sum of $n$ i.i.d.\ Bernoulli random variables with parameter $\alpha$, then
\begin{eqnarray}
|p_n^Z(x) - p_n^Z(x+z_0)| &=& \left|\sum_{j=1}^n \P(M_n=j) \sum_{z\in\Z^2} \left(p^{(1)}_{j}(z) - p^{(1)}_j(z+z_0)\right)
    p^{(2)}_{n-j}(x-z)\right| \nonumber\\
&\leq& \sum_{j=1}^n \P(M_n=j) \sum_{z\in\Z^2} \left|p^{(1)}_j(z)-p^{(1)}_j(z+z_0)\right| p^{(2)}_{n-j}(x-z). \label{A94}
\end{eqnarray}
Since $p^{(1)}$ has finite range, by (\ref{lclt}), it is easy to check that if $e$ is a unit vector in $\Z^2$ such 
that $x\cdot Q_1^{-1}x\leq (x+e)\cdot Q_1^{-1}(x+e)$, then
for any integer $k\geq 3$, we have
\be\label{diffest}
|p^{(1)}_n(x)-p^{(1)}_n(x+e)| \leq \frac{c}{n^{\frac{3}{2}}}\left[1+\left(\frac{\Vert x\Vert}{\sqrt n}\right)^{k+1}\right]
e^{-\frac{x\cdot Q_1^{-1}x}{2n}}+ o\big(n^{-\frac{k}{2}}\big) 
\ee
uniformly in $x$ and $n$. This bound and the local central limit theorem applied to $p^{(2)}$ are all we need to bound
(\ref{A94}) and establish (\ref{A888}).

Let $R=\max_{x\in\Z^2} \{\Vert x\Vert : p^{(1)}(x)>0\}$. Applying (\ref{diffest}) with $k=5$ then gives
\begin{eqnarray}
\sum_{x\in\Z^2} |p^{(1)}_n(x)-p^{(1)}_n(x+e)| &=&\sum_{\Vert x\Vert \leq Rn+1} |p^{(1)}_n(x)-p^{(1)}_n(x+e)| \nonumber\\
&\leq& \frac{2c}{\sqrt n} \sum_{\Vert x\Vert \leq Rn+1} \frac{1}{n} \left[1+ \left(\frac{\Vert x\Vert}{\sqrt
      n}\right)^6\right] e^{-\frac{x\cdot Q_1^{-1}x}{2n}} + \frac{C}{\sqrt n} \nonumber \\
&\leq& \frac{C}{\sqrt n} \label{A96},
\end{eqnarray}
where on the second line, the factor $2$ takes care of the possibility that $(x+e)\cdot Q_1^{-1}(x+e)< x\cdot Q_1^{-1}x$,
$C$ is uniform in $n\in\N$, and for the last inequality we used the Riemann sum approximation. By the triangle
inequality,
\be\label{diffest2}
\sum_{x\in\Z^2} |p_n^{(1)}(x) - p_n^{(1)}(x+z_0)| \leq \frac{C \Vert z_0\Vert}{\sqrt n}
\ee
with $C$ uniform in $z_0\in\Z^2$ and $n\in\N$. Using the decomposition (\ref{ppcouple}), it is easy to check that 
(\ref{diffest2}) in fact holds for all irreducible aperiodic random walks on $\Z^d$ (see Proposition 2.4.2 in
\cite{LL08}).

By the symmetry of (\ref{A888}) in $x$ and $x+z_0$, we may assume without loss of generality that
\be\label{asymmetry}
x\cdot Q_1^{-1}x \leq (x+z_0)\cdot Q_1^{-1}(x+z_0).
\ee 

To bound $\sum_n |p^Z_n(x) - p^Z_n(x+z_0)|$, we separate the sum into three regimes: 
\begin{itemize}
\item[(1)] $n\geq (1+\Vert x\Vert)^2$; 
\item[(2)] $1\leq n< c \frac{(1+\Vert x\Vert)^2}{\log (2+\Vert x\Vert)}$ for some $c>0$ sufficiently small; 
\item[(3)] $c\frac{(1+\Vert x\Vert)^2}{\log (2+\Vert x\Vert)} \leq n < (1+\Vert x\Vert)^2$.
\end{itemize}

For the regime $n\geq (1+\Vert x\Vert)^2$, by (\ref{A94}), (\ref{diffest2}) and the local central limit theorem for $p^{(2)}$, 
\begin{eqnarray}
|p^Z_n(x)-p^Z_n(x+z_0)| &\leq& \sum_{j=1}^n \P(M_n=j) \frac{C\Vert z_0\Vert}{\sqrt j}\frac{C}{1+n-j} \label{A98} \\
&\leq& \Vert z_0\Vert \left(C\P\Big(|M_n-\alpha n|\geq \alpha n/2\Big) + \frac{C}{n^{3/2}}\right) \nonumber \\
&\leq& C \frac{\Vert z_0\Vert }{n^{3/2}}, \nonumber
\end{eqnarray}
where we used elementary large deviation estimates for $M_n/n$. Therefore
\be\label{A99}
\sum_{n=(1+\Vert x\Vert)^2}^\infty |p^Z_n(x)-p^Z_n(x+z_0)| \leq \sum_{n=(1+\Vert x\Vert)^2}^\infty C \frac{\Vert
  z_0\Vert}{n^{3/2}} \leq C \frac{\Vert z_0\Vert}{1+\Vert x\Vert}
\ee
for some $C$ uniform in $x, z_0\in\Z^2$.

Now let $1\leq n \leq c(1+\Vert x\Vert)^2/\log (2+\Vert x\Vert)$ with $c>0$ to be chosen later. By our assumption 
(\ref{asymmetry}), we have $\Vert x+z_0\Vert \geq 2\epsilon \Vert x\Vert$ for some $\epsilon\in (0,1/2)$ depending only
on the smallest and largest eigenvalues of $Q_1$. Since $p^{(1)}$ has mean zero and finite range, by Hoeffding's 
concentration inequality \cite{H63} for martingales with bounded increments, uniformly for all
$1\leq j\leq c(1+\Vert x\Vert)^2/\log(2+\Vert x\Vert)$, we have
\be\label{A100}
\sum_{\Vert z\Vert \geq  \epsilon \Vert x\Vert } p^{(1)}_j(z) \leq C_1 e^{-C_2 \frac{\epsilon^2 \Vert x\Vert^2}{j}}
\leq C_1 e^{-C_2 \frac{\epsilon^2 \Vert x\Vert^2\log(2+\Vert x\Vert)}{c(1+\Vert x\Vert)^2}} \leq \frac{C}{(1+\Vert x\Vert)^3}
\ee
provided we choose $c<C_2\epsilon^2/3$. By (\ref{A94}),
\begin{eqnarray}
|p_n^Z(x) - p^Z_n(x+z_0)| &\leq& \sum_{j=1}^n P(M_n=j) \sum_{\Vert z-x\Vert \leq \epsilon \Vert x\Vert}
\left|p^{(1)}_j(z)-p^{(1)}_j(z+z_0)\right| p^{(2)}_{n-j}(x-z) \label{A101} \\
&& \quad + \sum_{j=1}^n \P(M_n=j) \sum_{\Vert z-x\Vert > \epsilon \Vert x\Vert} \left|p^{(1)}_j(z)-p^{(1)}_j(z+z_0)\right|
p^{(2)}_{n-j}(x-z) . \nonumber
\end{eqnarray}
Since we have assumed $\Vert x+z_0\Vert \geq 2\epsilon \Vert x\Vert$ for some $\epsilon\in (0,1/2)$, 
$\Vert z-x\Vert \leq \epsilon \Vert x\Vert$ implies $\Vert z\Vert \geq (1-\epsilon)\Vert x\Vert\geq \epsilon \Vert
x\Vert$ and $\Vert z+z_0\Vert =\Vert (x+z_0)+(z-x)\Vert \geq \epsilon \Vert x\Vert$. Therefore by (\ref{A100}), the 
first sum in (\ref{A101}) is bounded by $\frac{2C}{(1+\Vert x\Vert)^3}$. On the other hand, we have the following 
version of local central limit theorem for $p^{(2)}$ (see Section 7, P10 of \cite{S76}),
\be\label{A102}
p^{(2)}_j(y) = \frac{1}{2\pi j \sqrt{\det Q_2}} \left(e^{-\frac{y\cdot
      Q_2^{-1}y}{2j}}+o(1)\frac{j}{(1+\Vert y\Vert)^2}\right) \leq \frac{C}{(1+\Vert y\Vert)^2},
\ee  
where $C$ is uniform in $j$ and $y$. Combined with (\ref{diffest2}) and large deviation estimates for $M_n/n$,
this implies that the second sum in (\ref{A101}) is
bounded by $\frac{C \Vert z_0\Vert}{\sqrt{n}(1+\Vert x\Vert)^2}$ for some constant $C$ depending only $p^Z, p^{(1)}$ and
$p^{(2)}$. Therefore
\be\label{A103}
\sum_{n=1}^{c\frac{(1+\Vert x\Vert)^2}{\log(2+\Vert x\Vert)}}\!\!\!\!\! |p^Z_n(x) -p^Z_n(x+z_0)| \leq  \frac{2cC(1+\Vert x\Vert)^2}{(1+\Vert
  x\Vert)^3\log(2+\Vert x\Vert)} + \frac{C \Vert z_0\Vert}{(1+\Vert x\Vert)^2} \sum_{n=1}^{c\frac{(1+\Vert x\Vert)^2}{\log(2+\Vert
    x\Vert)}}\!\!\!\!\! \frac{1}{\sqrt n} \leq \frac{C\Vert z_0\Vert}{1+\Vert x\Vert}.
\ee

Finally, we treat the regime $c(1+\Vert x\Vert)^2/\log(2+\Vert x\Vert) \leq n \leq (1+\Vert x\Vert)^2$. By large deviation
estimates for $M_n/n$, it is easy to verify that 
\be\label{A107}
\sum_{n=c\frac{(1+\Vert x\Vert)^2}{\log(2+\Vert x\Vert)}}^{(1+\Vert x\Vert)^2} \sum_{1\leq j\leq n, \atop|j-\alpha
  n|\geq \alpha n/2}\!\!\!\! \P(M_n=j) \sum_{z\in\Z^2}\left|p^{(1)}_j(z)-p^{(1)}_j(z+z_0)\right| p^{(2)}_{n-j}(x-z) \leq \frac{C}{1+\Vert x\Vert}.
\ee
So we focus on $\alpha n/2 \leq j\leq 3\alpha n/2$ in (\ref{A94}). 

By the local central limit theorem for $p^{(2)}$, we have $p^{(2)}_i(y) \leq \frac{C}{i}$ uniformly for all 
$y\in\Z^2$ and $i\in\N$. Combined with (\ref{A102}), this implies that
\be\label{lcltcomb}
p^{(2)}_i(y) \leq \frac{C}{i \vee (1+\Vert y\Vert)^2}
\ee
for some $C$ uniformly in $y\in\Z^2$ and $i\in\N$. Therefore for all $\alpha n/2\leq j\leq 3\alpha n /2$ and $x,z\in\Z^2$,
we have
\be
p^{(2)}_{n-j}(x-z) \leq \frac{C}{n \vee (1+\Vert x-z\Vert)^2}.
\ee

If $v_0=0, v_1, \cdots, v_L=z_0$ is a nearest neighbor path in $\Z^2$ from $0$ to $z_0$, then by similar computations
as those leading to (\ref{A96}) except we now apply (\ref{diffest}) with $k=8$, we get
\begin{eqnarray}
&& \sum_{z\in\Z^2} \left|p^{(1)}_j(z)-p^{(1)}_j(z+z_0)\right| p^{(2)}_{n-j}(x-z) \nonumber \\
&\leq& \sum_{r=1}^L \sum_{z\in\Z^2} \left|p^{(1)}_j(z+v_{r-1})-p^{(1)}_j(z+v_r)\right| \frac{C}{n \vee (1+\Vert
  x-z\Vert)^2} \nonumber\\
&\leq& \sum_{r=1}^L \left(\sum_{z\in\Z^2} \frac{C}{j^{3/2}}\left[1+\left(\frac{\Vert z+v_r\Vert}{\sqrt
        j}\right)^9\right]\frac{e^{-\frac{(z+v_r)\cdot Q_1^{-1}(z+v_r)}{2j}}}{n \vee (1+\Vert x-z\Vert)^2} +
  \frac{C}{j^2}\right) \nonumber\\
&\leq& \sum_{r=1}^L \sum_{\tilde z\in\Z^2} \frac{C}{n^{3/2}}\left[1+\left(\frac{\Vert \tilde z\Vert}{\sqrt
        n}\right)^9\right]\frac{e^{-\frac{\tilde z\cdot Q_1^{-1}\tilde z}{2n}}}{n \vee (1+\Vert x+v_r-\tilde z\Vert)^2} +
  \frac{CL}{n^2}. \label{A108} 
\end{eqnarray}
For any $y\in\Z^2$, we have
\begin{eqnarray}
\!\!\!\!\!\!\!\!\!\!\!\!\!\!\!
&& \sum_{\tilde z\in\Z^2} \frac{C}{n^{3/2}}\left[1+\left(\frac{\Vert \tilde z\Vert}{\sqrt
        n}\right)^9\right]\frac{e^{-\frac{\tilde z\cdot Q_1^{-1}\tilde z}{2n}}}{n \vee (1+\Vert y-\tilde z\Vert)^2}
  \nonumber \\
\!\!\!\!\!\!\!\!\!\!\!\!\!\!\!
&\leq& \sum_{\Vert\tilde z\Vert \geq \frac{\Vert y\Vert}{2}} \frac{C}{n^{5/2}} \left[1+\left(\frac{\Vert \tilde z\Vert}{\sqrt
        n}\right)^9\right]e^{-\frac{\tilde z\cdot Q_1^{-1}\tilde z}{2n}} + \frac{C}{\big(1+\frac{\Vert y\Vert}{2}\big)^2} 
\sum_{\Vert\tilde
    z\Vert <\frac{\Vert y\Vert}{2}} \frac{1}{n^{3/2}}
\left[1+\left(\frac{\Vert \tilde z\Vert}{\sqrt n}\right)^9\right]e^{-\frac{\tilde z\cdot Q_1^{-1}\tilde z}{2n}}. \label{A109}
\end{eqnarray}
Note that by Riemann sum approximation, 
$$
\sum_{\tilde z\in\Z^2}\frac{1}{n}\left[1+\left(\frac{\Vert \tilde z\Vert}{\sqrt
      n}\right)^9\right]e^{-\frac{\tilde z\cdot Q_1^{-1}\tilde z}{4n}} \underset{n\to\infty}{\longrightarrow}
\int_{\R^2} (1+\Vert w\Vert^9)e^{-\frac{w\cdot Q_1^{-1} w}{4}}dw <\infty.
$$
It is then easy to see that in (\ref{A109}),  the first sum is bounded by $\frac{C_1}{n^{3/2}} e^{-C_2\frac{\Vert
y\Vert^2}{n}}$ and the second sum is bounded by $\frac{C_3}{\sqrt{n}(1+\Vert y\Vert)^2}$, where $C_1,
C_2$ and $C_3$ are uniform in $y\in\Z^2$ and $n\in\N$. Hence
\be\label{A110}
\!\!\!\!\sum_{\tilde z\in\Z^2}\! \frac{C}{n^{3/2}}\left[1+\left(\frac{\Vert \tilde z\Vert}{\sqrt
        n}\right)^9\right]\frac{e^{-\frac{\tilde z\cdot Q_1^{-1}\tilde z}{2n}}}{n \vee (1+\Vert y-\tilde z\Vert)^2}
\leq \frac{C_1}{n^{3/2}} e^{-C_2\frac{\Vert y\Vert^2}{n}} + \frac{C_3}{\sqrt{n}(1+\Vert y\Vert)^2}
\leq \frac{C}{\sqrt{n} (1+\Vert y\Vert)^2}.
\ee
By our assumption (\ref{asymmetry}), which implies $\Vert x+z_0\Vert \geq 2\epsilon \Vert x\Vert$ for some $\epsilon \in
(0,1/2)$ depending only on $Q_1$, we can choose the nearest neighbor path $v_0=0, v_1, \cdots, v_L=z_0$ such that 
$L\leq C\Vert z_0\Vert$ for some $C$ independent of $x$ and $z_0$, and $\Vert x+v_r\Vert \geq \epsilon \Vert x\Vert$ for
all $0\leq r\leq L$. For such a path, we can substitute the bound (\ref{A110}) into (\ref{A108}) to obtain
\be
\sum_{|j-\alpha n|< \alpha n/2}\P(M_n=j) \sum_{z\in\Z^2} \left|p^{(1)}_j(z)-p^{(1)}_j(z+z_0)\right| p^{(2)}_{n-j}(x-z) \leq \frac{C\Vert
  z_0\Vert}{\sqrt{n}(1+\epsilon\Vert x\Vert)^2} + \frac{C\Vert z_0\Vert}{n^2}.
\ee 
Combined with (\ref{A107}), we see that
\be
\sum_{n=c\frac{(1+\Vert x\Vert)^2}{\log(2+\Vert x\Vert)}}^{(1+\Vert x\Vert)^2} \left| p^Z_n(x)-p^Z_n(x+z_0)\right| \leq
\frac{C\Vert z_0\Vert}{1+\Vert x\Vert}
\ee
with $C$ uniform in $x, z_0\in\Z^2$. Together with (\ref{A99}) and (\ref{A103}), this proves (\ref{A888}).
\qed
\bigskip

To prove Lemma \ref{L:rwconv}, we will use the so-called rearrangement inequality. For much deeper results on
rearrangement inequalities than the one we use here, see Chapter 3 of Lieb and Loss \cite{LL01}.

\bl\label{L:rearrange} Let $(a_n)_{n\in\N}$ be a non-negative non-increasing sequence, and let $(b_n)_{n\in\N}$
and $(c_n)_{n\in\N}$ be two non-negative sequences. If $c$ majorizes $b$ in the sense that $\sum_{i=1}^n b_i \leq
\sum_{i=1}^n c_i$ for all $n\in\N$,
then 
\be\label{rearrange1}
\sum_{n=1}^\infty a_n b_n \leq \sum_{n=1}^\infty a_n c_n.
\ee
In particular, if there exists a bijection $\sigma :\N\to\N$ such that $(b_{\sigma(n)})_{n\in\N}$ becomes a
non-increasing sequence, then
\be\label{rearrange2}
\sum_{n=1}^\infty a_n b_n \leq \sum_{n=1}^\infty a_n b_{\sigma(n)}.
\ee
\el
The proof of Lemma \ref{L:rearrange} is elementary, so we omit it. The majorization condition can be interpreted
as a stochastic domination condition for the positive measures on $\N$ defined by $b$ and $c$.
\bigskip

\noindent
{\bf Proof of Lemma \ref{L:rwconv}.} We only treat the discrete time 
aperiodic case. The discrete time periodic case and the continuous time case 
are similar. Let $(x_n)_{n\in\N}$ be an ordering of $\Z^2$ in increasing 
Euclidean norm. 
Clearly the sequence $\frac{1}{(1+\Vert x_n\Vert)^q}$, $n\in\N$, is non-increasing. Let $(y_n)_{n\in\N}$ be
an ordering of $\Z^2$ such that $\P(X_i=y_n)$ becomes a non-increasing sequence. Then by the rearrangement
inequality (\ref{rearrange2}), 
\be\label{A88}
\sum_{x\in\Z^2} \P(X_i=x) \frac{1}{(1+\Vert x-v\Vert)^q} \leq \sum_{n=1}^\infty \P(X_i=y_n)\frac{1}{(1+\Vert x_n\Vert)^q}.
\ee
Let $Q$ denote the covariance matrix of $X$. By the local central limit theorem, 
$\P(X_i=x) = \frac{e^{-\frac{\langle x, Q^{-1}x\rangle}{2i}}+o(1)}{2\pi i \sqrt{\det Q}}$ uniformly in $x$. Since
$Q^{-1}$ is positive definite, we can choose $C$ and $\alpha$ independent of $i$, such that $(\P(X_i=y_n))_{n\in\N}$
is majorized by (as defined in Lemma \ref{L:rearrange}) the sequence $(b_n)_{n\in\N}$ with 
$b_n = \frac{C e^{-\frac{\alpha\Vert x_n\Vert^2}{i} }}{i}$ 
when $\Vert x_n\Vert \leq \sqrt{i}$, and $b_n=0$ when $\Vert x_n\Vert >\sqrt{i}$. Then by (\ref{rearrange1}),
\be\label{A89}
\sum_{n=1}^\infty \P(X_i=y_n)\frac{1}{(1+\Vert x_n\Vert)^q} \leq \sum_{\Vert x\Vert \leq \sqrt{i}} C\frac{
  e^{-\frac{\alpha\Vert x\Vert^2}{i}}}{i (1+\Vert x\Vert)^q}.
\ee
By Riemann sum approximation,
\be
\lim_{i\to\infty} i^{q/2} \sum_{\Vert x\Vert \leq \sqrt{i}} C\frac{e^{-\frac{\alpha\Vert x\Vert^2}{i}}}{i (1+\Vert x\Vert)^q}
=\lim_{i\to\infty} \sum_{\Vert x\Vert \leq \sqrt{i}} C\frac{e^{-\frac{\alpha\Vert
      x\Vert^2}{i}}}{(\frac{1}{\sqrt{i}}+\frac{\Vert x\Vert}{\sqrt i})^q}\ \frac{1}{i} = C\int_{\Vert w\Vert\leq
  1}\frac{e^{-\alpha\Vert w\Vert^2}}{\Vert w\Vert^q}dw,
\ee
which is finite if $q<2$. In view of (\ref{A88}) and (\ref{A89}), this implies (\ref{A72}). 
\qed

\section{Proof of Theorem \ref{T:2rw}}\label{S3}
Since the proof of Theorem \ref{T:explaw} is rather simple, we include it here for completeness.

\medskip
\noindent
{\bf Proof of Theorem \ref{T:explaw}.} We give the proof for the continuous time random walk case. The discrete 
time case can be treated similarly, or it can be deduced from the continuous time case by a change of time argument.
Let $p_t(x)=\P(X_t=x|X_0=0)$. Note that for each $k\in\N$,
\begin{eqnarray}
\E\left[L_t^k\right] &=& \int_0^t ds_1\cdots \int_0^t ds_k \P(X_{s_1}=\cdots X_{s_k}=0) \nonumber\\
&=& k! \int_0^t ds_1 p_{s_1}(0) \int_0^{t-s_1}ds_2 p_{s_2}(0) \cdots \int_0^{t-\sum_{i=1}^{k-1}s_i}
p_{s_k}(0)ds_k. \nonumber
\end{eqnarray}
Clearly 
\be
\left(\int_0^{t/k} p_s(0)ds\right)^k\ \leq\ \frac{\E[L_t^k]}{k!}\ \leq\ \left(\int_0^t p_s(0)ds\right)^k.
\ee
By the local central limit theorem, $p_t(0) \sim \frac{1}{2\pi rt \sqrt{\det Q}}$, where we write $f(t)\sim g(t)$
if $\lim_{t\to\infty} f/g =1$. Therefore $\int_0^{ct} p_s(0)ds \sim \frac{\log t}{2\pi r\sqrt{\det Q}}$ for any $c>0$, 
and 
\be\label{kthmoment}
\lim_{t\to\infty} \E\left[\Big(\frac{2\pi r\sqrt{\det Q}L_t}{\log t}\Big)^k\right] = k!\ .
\ee
Since $k!$ is the $k$-th moment of a mean 1 exponential random variable and is distribution determining, the desired 
convergence in distribution follows by the method of moments.
\qed

\bigskip
\noindent
{\bf Proof of Theorem \ref{T:2rw}.}
For simplicity, we write $L_t$ for $L_t(X,Y)$ from now on. We divide the proof into three parts. First 
we treat the discrete time case and show that for each $k\in\N$ and $\epsilon>0$,
\be\label{varbound}
{\rm Var}\big(\E^X_0[L_n^k|Y]\big) = \E^Y_0\Big[\Big(\E^X_0[L_n^k|Y] - E^{X,Y}_{0,0}[L_n^k]\Big)^2\Big] = o(\log^{2k-1+\epsilon} n),
\ee
which implies a weak law of large numbers for $\E^X_0[L_n^k|Y]/\log^k n$ as $n\to\infty$. We then show how to
adapt the argument to the continuous time case. Lastly, we show that our variance bounds in fact imply a strong law of
large numbers for $\E^X_0[L_n^k|Y]/\log^k n$. The claimed almost sure convergence in distribution for
$\frac{2\pi(\kappa+\rho)\sqrt{\det Q} L_n(X,Y)}{\log n}$ then follows by the method of moments.

\bigskip
\noindent
{\bf Variance bound for discrete time random walks.} Let $F_n$ denote the sigma-field generated by $(Y_i)_{0\leq i\leq n}$. By
the martingale decomposition, for any $f_n(Y)$ measurable w.r.t.\ $F_n$, we have
\be\label{martdec}
\E^Y_0[(f_n-E^Y_0[f_n])^2] = \sum_{i=1}^n \E^Y_0\left[\big(\E^Y_0[f_n|F_i]-\E^Y_0[f_n|F_{i-1}]\big)^2\right].
\ee
We now estimate the $i$-th term in the summation. Note that $(\E^Y_0[f_n|F_i]-\E^Y_0[f_n|F_{i-1}])^2$ depends only on
$(Y_j)_{1\leq j\leq i}$. Let us first integrate out the last jump $Y_i-Y_{i-1}$. By the standard trick that 
$\E[(Z-\E[Z])^2]= \frac{1}{2}\E[(Z-Z')^2]$ where $Z'$ is an independent copy of $Z$, we have
\begin{eqnarray}
&& \E^Y_0\left[\big(\E^Y_0[f_n|F_i]-\E^Y_0[f_n|F_{i-1}]\big)^2\,\Big|\, F_{i-1}\right] \nonumber \\
&=& \frac{1}{2}\E^{\Delta, \Delta'}\left[\Big(\E^Y_0[f_n | F_{i-1}, Y_i-Y_{i-1}=\Delta]-\E^Y_0[f_n | F_{i-1},
  Y_i-Y_{i-1}=\Delta']\Big)^2\right], \label{A59}
\end{eqnarray}
where $\Delta$ and $\Delta'$ are independent copies of the increment of $Y$ in one step, and hence
\begin{eqnarray}
\!\!\!\!\!\!&& \E^Y_0\left[\big(\E^Y_0[f_n|F_i]-\E^Y_0[f_n|F_{i-1}]\big)^2\right] \nonumber \\
\!\!\!\!\!\!&=& \frac{1}{2}\E^{\Delta, \Delta'}\left[ \E^Y_0\left[\Big(\E^Y_0[f_n | F_{i-1}, Y_i-Y_{i-1}=\Delta]-\E^Y_0[f_n | F_{i-1},
    Y_i-Y_{i-1}=\Delta']\Big)^2 \Big| \Delta, \Delta'\right]\right]. \label{2copy}
\end{eqnarray}

We now specialize to the case $f_n(Y) = \E^X_0[L^k_n(X,Y)|Y]$ for some fixed $k\in\N$. Write $L_n= L_{[0,i-1]} + L_{[i,n]}$ where 
$L_{[a,b]} = \sum_{a\leq j\leq b}\delta_0(X_j-Y_j)$. Then 
\be \label{Ln}
L_n^k = L_{[0,i-1]}^k + \sum_{m=1}^k \binom{k}{m} L_{[0,i-1]}^{k-m} L_{[i,n]}^m.
\ee
Write $\Delta$, resp.\ $\Delta'$, as a shorthand for the conditioning $Y_i-Y_{i-1}=\Delta$, resp.\ $\Delta'$, and let
$p_i^X(\cdot)$ denote the $i$-step transition probability kernel for $X$. Then
\begin{eqnarray}
&& \E^Y_0[f_n | F_{i-1}, \Delta]-\E^Y_0[f_n | F_{i-1}, \Delta'] \nonumber \\
&=& \sum_{m=1}^k \binom{k}{m} \Big(\E^{X,Y}_{0,0}\big[
L_{[0,i-1]}^{k-m} L_{[i,n]}^m \,\big|\, F_{i-1}, \Delta\big] - \E^{X,Y}_{0,0}\big[L_{[0,i-1]}^{k-m} L_{[i,n]}^m \,\big|\, F_{i-1},\Delta'\big]\Big) \nonumber \\
&=& \sum_{m=1}^k \binom{k}{m} \sum_{x\in\Z^2} p^X_i(x)\ \E^X_0\big[L_{[0,i-1]}^{k-m} \,\big|\, F_{i-1}, X_i=x\big]  \label{A62}\\
&& \qquad  
\times \Big(\E^{X,Y}_{0,0}\big[L_{[i,n]}^m \,\big|\, F_{i-1}, \Delta, X_i=x\big] - \E^{X,Y}_{0,0}\big[L_{[i,n]}^m \,\big|\,F_{i-1}, \Delta',
X_i=x\big]\Big). \nonumber
\end{eqnarray}
If we denote $Y_{i-1}=y$, and denote $Z= X-Y$, then we have 
\begin{eqnarray}
&& \E^Y_0[f_n | F_{i-1}, \Delta]-\E^Y_0[f_n | F_{i-1}, \Delta'] \label{A63} \\
&=& \sum_{m=1}^k \binom{k}{m} \sum_{x\in\Z^2} p_i^X(x)\ \E^X_0[L_{[0,i-1]}^{k-m} |F_{i-1}, X_i=x]\ \Big(\E^Z_{x-y-\Delta}[L_{n-i}^m(Z)]
-\E^Z_{x-y-\Delta'}[L_{n-i}^m(Z)]\Big), \nonumber
\end{eqnarray}
where $L_n(Z)=\sum_{j=0}^n\delta_0(Z_j)$. It is easy to see that
\begin{eqnarray}
&& \Big(\E^Y_0[f_n | F_{i-1}, \Delta]-\E^Y_0[f_n | F_{i-1}, \Delta']\Big)^2 \label{A67} \\
&\leq& C_{k}\sum_{m=1}^k \left(\sum_{x\in\Z^2} p_i^X(x)\ \E^X_0[L_{[0,i-1]}^{k-m} |F_{i-1}, X_i=x]\ \Big|\E^Z_{x-y-\Delta}[L_{n-i}^m(Z)]
-\E^Z_{x-y-\Delta'}[L_{n-i}^m(Z)]\Big|\right)^2. \nonumber
\end{eqnarray}
Here and as well as in what follows, $C_k$ always denotes a generic constant depending only on $k$ and the transition
kernels of $X$ and $Y$, whose precise value may change from line to line.

By expanding $L_n^m(Z) = (\sum_{0\leq j\leq n}\delta_0(Z_j))^m$, we have
\begin{eqnarray}
&& \Big|\E^Z_{x-y-\Delta}[L_{n-i}^m(Z)] -\E^Z_{x-y-\Delta'}[L_{n-i}^m(Z)]\Big| \nonumber\\
&\leq& \sum_{0\leq j_1, \cdots, j_m\leq n-i} \Big| \P^Z_{x-y-\Delta}(Z_{j_1}= \cdots
  Z_{j_m}=0)-\P^Z_{x-y-\Delta'}(Z_{j_1}=\cdots =Z_{j_m}=0)\Big| \nonumber \\
&\leq& m! \sum_{0\leq j_1\leq j_2\cdots\leq j_{m}\leq n} \big|p^Z_{j_1}(x-y-\Delta)-p^Z_{j_1}(x-y-\Delta')\big|
p^Z_{j_2-j_1}(0)\cdots p^Z_{j_m-j_{m-1}}(0) \nonumber\\
&\leq& m! \sum_{0\leq j_1< \infty} \big|p^Z_{j_1}(x-y-\Delta)-p^Z_{j_1}(x-y-\Delta')\big| \left(\sum_{0\leq j_2\leq
    n} p^Z_{j_2}(0)\right)^{m-1} \nonumber \\
&\leq& C_k (\log n)^{m-1} \sum_{0\leq j_1 < \infty} \big|p^Z_{j_1}(x-y-\Delta)-p^Z_{j_1}(x-y-\Delta')\big| \nonumber\\
&\leq& C_k (\log n)^{m-1} \Vert \Delta-\Delta'\Vert \left(\frac{1}{1+\Vert x-y-\Delta\Vert}+\frac{1}{1+\Vert
    x-y-\Delta'\Vert}\right), \label{A688}
\end{eqnarray}
where in the last two inequalities, we used the local central limit theorem which implies that $p^Z_n(0) \leq cn^{-1}$ for some $c>0$,
and we applied Lemma \ref{L:gradpot}. Substituting (\ref{A688}) into (\ref{A67}), we get
\begin{eqnarray}
&& \Big(\E^Y_0[f_n | F_{i-1}, \Delta]-\E^Y_0[f_n | F_{i-1}, \Delta']\Big)^2 \nonumber \\
&\leq& C_k \sum_{m=1}^k (\log n)^{2(m-1)}\Vert \Delta-\Delta'\Vert^2  \left(\sum_{v=\Delta,\Delta'} \sum_{x\in\Z^2}
  p_i^X(x) \frac{\E^X_0[L_{[0,i-1]}^{k-m} |F_{i-1}, X_i=x]}{1+\Vert x-y-v\Vert} \right)^2 \nonumber \\
&\leq& C_k \sum_{m=1}^k (\log n)^{2(m-1)}\Vert \Delta-\Delta'\Vert^2  \sum_{v=\Delta,\Delta'} \left(\sum_{x\in\Z^2}
  p_i^X(x) \frac{\E^X_0[L_{[0,i-1]}^{k-m} |F_{i-1}, X_i=x]}{1+\Vert x-y-v\Vert} \right)^2. \label{A71}
\end{eqnarray}

Let $q\in (1,2)$ and $p\in(2,\infty)$ with $\frac{1}{p}+\frac{1}{q}=1$. We now apply H\"older's inequality and Lemma
\ref{L:rwconv} to obtain

\begin{eqnarray}
&& \sum_{x\in\Z^2} p_i^X(x) \frac{\E^X_0[L_{[0,i-1]}^{k-m} |F_{i-1}, X_i=x]}{1+\Vert x-y-v\Vert} \nonumber \\
&\leq& \left(\sum_{x\in\Z^2} p_i^X(x) \frac{1}{(1+\Vert
    x-y-v\Vert)^q}\right)^{\frac{1}{q}}\Big(\sum_{x\in\Z^2} p_i^X(x) \left(\E^X_0[L_{[0,i-1]}^{k-m} |F_{i-1},
  X_i=x]\right)^p\Big)^{\frac{1}{p}} \nonumber \\
&\leq& \frac{C}{i^{\frac{1}{2}}} \Big(1+\sum_{x\in\Z^2} p_i^X(x) \left(\E^X_0[L_{[0,i-1]}^{k-m} |F_{i-1},
  X_i=x]\right)^p\Big)^{\frac{1}{2}}. \label{A73}
\end{eqnarray}
Substituting (\ref{A73}) into (\ref{A71}) then gives
\begin{eqnarray}
\!\!\!\!\!&& \Big(\E^Y_0[f_n | F_{i-1}, \Delta]-\E^Y_0[f_n | F_{i-1}, \Delta']\Big)^2 \nonumber \\
\!\!\!\!\!&\leq& C_k \sum_{m=1}^k (\log n)^{2(m-1)} \frac{\Vert\Delta-\Delta'\Vert^2}{i}\Big(1+\sum_{x\in\Z^2}
p_i^X(x) \left(\E^X_0[L_{[0,i-1]}^{k-m} |F_{i-1},
  X_i=x]\right)^p\Big).
\end{eqnarray}
The important point here is that we obtain a factor of $\frac{1}{i}$. We can now substitue this estimate into
(\ref{2copy}) with $f_n=\E^X_0[L_n^k(X,Y) |Y]$ to get
\begin{eqnarray}
\!\!\!\!\!\! &&\!\!\!\!\! 
\E^Y_0\left[\big(\E^Y_0[f_n|F_i]-\E^Y_0[f_n|F_{i-1}]\big)^2\right] \nonumber \\
\!\!\!\!\!\!\!\!\!\! &\leq&\!\!\!\!
\frac{C_k \E^{\Delta,\Delta'}[\Vert\Delta-\Delta'\Vert^2]}{i} \sum_{m=1}^k (\log
n)^{2(m-1)}\Big(1+\sum_{x\in\Z^2} p_i^X(x)\, \E^Y_0\left[\big(\E^X_0[L_{[0,i-1]}^{k-m} |F_{i-1},X_i=x]\big)^p\right]
\Big). \label{A75} 
\end{eqnarray}
Since $p>1$, applying Minkowski's inequality (an integral version of the triangle inequality on $L_p$ space, 
see Section 2.4 of Lieb and Loss \cite{LL01})
\be\label{Minkowski}
\left(\int_\Omega\Big(\int_{\Gamma}|f(x,y)|\nu(dx)\Big)^p \mu(dy)\right)^{\frac{1}{p}}
\leq \int_\Gamma \Big(\int_\Omega |f(x,y)|^p \mu(dy)\Big)^{\frac{1}{p}} \nu(dx)
\ee
to the two fold expectation in (\ref{A75}) with $\E^Y_0[\cdot]$ playing the role of $\int\cdot \mu(dy)$ and
$\E^X_0[\cdot|X_i=x]$ playing the role of $\int \cdot \nu(dx)$, we get
\be
\E^Y_0\left[\Big(\E^X_0\big[L_{[0,i-1]}^{k-m} \big|F_{i-1},X_i=x\big]\Big)^p\right] \leq \E^X_0\left[ \E^Y_0\left[L_{[0,i-1]}^{(k-m)p}\Big|(X_j)_{0\leq
    j\leq i}\right]^{\frac{1}{p}} \Big|X_i=x\right]^p. \label{A76}
\ee
The advantage of estimating the RHS of (\ref{A76}) over the LHS is that, we have a good uniform bound on
$\E^Y_0\left[L_{[0,i-1]}^{(k-m)p}\big|(X_j)_{1\leq j\leq i}\right]$ with respect to $(X_j)_{1\leq j\leq i}$, which allows us
to circumvent the conditioning on $X_i=x$. More precisely, by the local central limit theorem for $Y$, we have $p^Y_n(y) \leq
\frac{C}{n}$ for some $C$ uniformly in $n\in\N$ and $y\in\Z^2$. Hence by H\"older's inequality and the same expansion of
$L_{[0,n]}^m$ as the one leading to (\ref{A688}), we get
\be\label{A722}
\E^Y_0\left[L_{[0,i-1]}^{(k-m)p}\big|(X_j)_{0\leq j\leq i}\right] \leq 
\E^Y_0\left[L_{[0,n]}^{\lceil (k-m)p\rceil}\big|(X_j)_{0\leq j\leq n}\right]^{\frac{(k-m)p}{\lceil (k-m)p\rceil}}
\leq C_{k,p} (\log n)^{(k-m)p}
\ee 
uniformly in $(X_j)_{0\leq j\leq n}$. Substituting this bound into (\ref{A76}), (\ref{A75}) and then into (\ref{martdec}), and
combining various constants together, we get for $f_n(Y) = \E^X_0[L_n^k|Y]$, 
\begin{eqnarray}
\E^Y_0[(f_n-E^Y_0[f_n])^2] \leq C \sum_{i=1}^n \frac{1}{i} \sum_{m=1}^k (\log n)^{2(m-1)+(k-m)p} 
\leq C \sum_{m=1}^k(\log n)^{2m-1+(k-m)p},
\end{eqnarray}
where $C$ depends only on $p, k, X$ and $Y$. Since $p>2$ can be chosen to be arbitrarily close to $2$, we see that
${\rm Var}(\E^X_0[L_n^k|Y]) = o(\log^{2k-1+\epsilon}n)$ for all $\epsilon>0$, which is what we set out to prove.

\bigskip
\noindent{\bf Variance bound for continuous time random walks.} We now adapt the above argument to continuous time
random walks, which is a bit more cumbersome. Without loss of generality, assume $t=n\in\N$. The martingale decomposition
(\ref{martdec}) is still valid. However, in (\ref{A59}) and (\ref{2copy}), instead of conditioning on
$Y_i-Y_{i-1}=\Delta$, resp.\ $\Delta'$, we need to condition on $(Y_{i-1+s}-Y_{i-1})_{s\in
  [0,1]}=(\Delta_s)_{s\in[0,1]}$, resp.\ $(\Delta'_s)_{s\in[0,1]}$. We also need to replace (\ref{Ln}) by 
\be
L_n^k=L_{[0,i-1]}^k + \sum_{m=1}^k \binom{k}{m} L_{[0, i-1]}^{k-m} L_{[i-1,n]}^k.
\ee
To extract $L_{[0,i-1]}^{k-m}$ as a common factor as in (\ref{A62}), we should now condition on $X_{i-1}=x$ rather than
on $X_i=x$. Writing simply $\Delta$ as a shorthand for the conditioning
$(Y_{i-1+s}-Y_{i-1})_{s\in[0,1]}=(\Delta_s)_{s\in[0,1]}$, and the same for $\Delta'$, (\ref{A62}) is now replaced by
\begin{eqnarray}
&& \E^Y_0[f_n | F_{i-1}, \Delta]-\E^Y_0[f_n | F_{i-1}, \Delta'] \nonumber \\
&=& \sum_{m=1}^k \binom{k}{m} \sum_{x\in\Z^2} p_{i-1}^X(x)\ \E^X_0[L_{[0,i-1]}^{k-m} |F_{i-1}, X_{i-1}=x] \nonumber \\
&& \qquad \times \Big(\E^{X,Y}_{0,0}[L_{[i-1,n]}^m |F_{i-1}, \Delta, X_{i-1}=x] - \E^{X,Y}_{0,0}[L_{[i-1,n]}^m |F_{i-1},
\Delta', X_{i-1}=x]\Big). \label{A622} 
\end{eqnarray}
In (\ref{A622}), we make the further expansion that 
\be \label{Ln2}
L_{[i-1,n]}^m = \sum_{l=0}^m \binom{m}{l} L_{[i-1,i]}^l L_{[i,n]}^{m-l}.
\ee 
The resulting expansion for (\ref{A622}) then consists of the following three types of terms:
\begin{eqnarray}
\Gamma^{\Delta,\Delta'}_{m,0} &:& \sum_{x\in\Z^2} p_{i-1}^X(x)\ \E^X_0[L_{[0,i-1]}^{k-m}\, \big|\, F_{i-1}, X_{i-1}=x] \nonumber \\
&&  \qquad \qquad \times \Big(\E^{X,Y}_{0,0}[L_{[i,n]}^m \, \big|\, F_{i-1}, \Delta, X_{i-1}=x] -
\E^{X,Y}_{0,0}[L_{[i,n]}^m \, \big|\, F_{i-1}, \Delta', X_{i-1}=x]\Big), \nonumber \\
\Gamma^{\Delta}_{m,l},\ l\geq 1 &:& \sum_{x\in\Z^2} p_{i-1}^X(x)\ \E^X_0[L_{[0,i-1]}^{k-m}\, \big|\, F_{i-1}, X_{i-1}=x] 
\ \E^{X,Y}_{0,0}\Big[L_{[i-1,i]}^l L_{[i,n]}^{m-l} \, \big|\, F_{i-1}, \Delta, X_{i-1}=x\Big], \nonumber \\
\Gamma^{\Delta'}_{m,l},\ l\geq 1 &:& \sum_{x\in\Z^2} p_{i-1}^X(x)\ \E^X_0[L_{[0,i-1]}^{k-m} \, \big|\, F_{i-1}, X_{i-1}=x] 
\ \E^{X,Y}_{0,0}\Big[L_{[i-1,i]}^l L_{[i,n]}^{m-l} \,\big|\, F_{i-1}, \Delta', X_{i-1}=x\Big]. \nonumber
\end{eqnarray}
Therefore
\begin{eqnarray}
&& \Big(\E^Y_0[f_n | F_{i-1}, \Delta]-\E^Y_0[f_n | F_{i-1}, \Delta']\Big)^2 \nonumber \\
&=&  \left(\sum_{m=1}^k \binom{k}{m}\Gamma^{\Delta, \Delta'}_{m,0}
  +\sum_{m=1}^k\sum_{l=1}^m\binom{k}{m}\binom{m}{l}\Gamma^\Delta_{m,l} - \sum_{m=1}^k\sum_{l=1}^m
  \binom{k}{m}\binom{m}{l}\Gamma^{\Delta'}_{m,l}\right)^2   \nonumber\\
&\leq& C_k \left(\sum_{m=1}^k \Big(\Gamma^{\Delta, \Delta'}_{m,0}\Big)^2 +\sum_{m=1}^k\sum_{l=1}^m
  \Big(\Gamma^\Delta_{m,l}\Big)^2 + \sum_{m=1}^k\sum_{l=1}^m \Big(\Gamma^{\Delta'}_{m,l}\Big)^2 \right),
\end{eqnarray}
where $C_k$ is a constant depending only on $k$. To bound the variance as in (\ref{A59}), we need to bound
$$
\E^{\Delta,\Delta'}\Big[\,\E^Y_0\big[\big(\Gamma_{m,0}^{\Delta,\Delta'}\big)^2\big|\Delta, \Delta'\big]\,\Big],\quad 
\E^{\Delta,\Delta'}\Big[\,\E^Y_0\big[\big(\Gamma_{m,l}^{\Delta}\big)^2 \big|\Delta, \Delta'\big]\,\Big] \quad \text{ and } \quad
\E^{\Delta,\Delta'}\Big[\,\E^Y_0\big[\,\big(\Gamma_{m,l}^{\Delta'}\big)^2 \big|\Delta, \Delta'\big]\,\Big].
$$
For terms involving $\Gamma^{\Delta, \Delta'}_{m,0}$, $1\leq m\leq k$, if we denote $Y_{i-1}=y$ and by further
conditioning on $X_i=x'$, we then have
\begin{eqnarray}
&& \Big|\E^{X,Y}_{0,0}[L_{[i,n]}^m \, \big|\, F_{i-1}, \Delta, X_{i-1}=x] - \E^{X,Y}_{0,0}[L_{[i,n]}^m \, \big|\, F_{i-1},
\Delta', X_{i-1}=x]\Big| \nonumber \\
&\leq& \sum_{x'\in \Z^2} p^X_1(x'-x) \Big|\E^{Z}_{x'-y-\Delta_1}\big[L_{n-i}^m(Z)\big]-\E^{Z}_{x'-y-\Delta'_1}\big[L_{n-i}^m(Z)\big]  \Big| 
\nonumber\\
&\leq& C_k (\log n)^{m-1} \Vert \Delta_1-\Delta'_1\Vert \sum_{v=\Delta_1, \Delta'_1} \sum_{x'\in \Z^2} \frac{p^X_1(x'-x)}{1+\Vert x'-y -v\Vert}, 
\end{eqnarray}
where $Z=X-Y$, and we followed the same computation as in (\ref{A688}). Since $X$ has finite second moments, by the 
Markov inequality, we have
\begin{eqnarray}
\sum_{x'\in\Z^2} \frac{p_1^X(x'-x)}{1+\Vert x'-y -v\Vert} &\leq& \P^X_0\Big(\Vert X_1\Vert \geq \frac{1+\Vert x-y -v\Vert}{2}\Big)
+ \sup_{\Vert x'-x\Vert < \frac{1+\Vert x-y -v\Vert}{2}} \frac{1}{1+\Vert x'-y -v\Vert} \nonumber \\
&\leq& \frac{C}{\big(\frac{1+\Vert x-y -v\Vert}{2}\big)^2} + \sup_{\Vert x'-x\Vert < \frac{1+\Vert x-y -v\Vert}{2}} 
\frac{1}{1+\Vert x-y -v\Vert -\Vert x'-x\Vert} \nonumber\\
&\leq& \frac{C'}{1+\Vert x-y-v\Vert},
\end{eqnarray}
where $C$ and $C'$ are constants depending only on $X$. This reduces the bound for $\big(\Gamma^{\Delta, \Delta'}_{m,0}\big)^2$ to
the same form as in (\ref{A71}). The calculations for the discrete time case then carry over, and we conclude
that the contribution of terms involving $\Gamma^{\Delta, \Delta'}_{m,0}$ to the variance of $\E^X_0[L_n^k|Y]$ is 
of order $o(\log^{2k-1+\epsilon} n)$ for all $\epsilon>0$. 

We now bound $\E^{\Delta,\Delta'}\Big[\,\E^Y_0\big[\big(\Gamma_{m,l}^{\Delta}\big)^2 \big|\Delta,
\Delta'\big]\,\Big]$, $1\leq l \leq m$. The case involving $\Gamma^{\Delta'}_{m,l}$ is identical. By first conditioning with
respect to $X_i=x'$ and then applying the local central limit theorem as in (\ref{A722}), and using the fact that
$L_{[i-1,i]}\leq 1$ and $l\geq 1$, we get
\be
\E^{X,Y}_{0,0}\Big[L_{[i-1,i]}^l L_{[i,n]}^{m-l} \Big| F_{i-1}, \Delta, X_{i-1}=x\Big] \leq C (\log n)^{m-l} \E^X_{x-Y_{i-1}}\big[L_{[0,1]}(X,\Delta)\big].
\ee
Hence
\begin{eqnarray}
\!\!\!\!\!\! && \label{A811}
(\log n)^{-2(m-l)}\, \E^{\Delta,\Delta'}\Big[\,\E^Y_0\big[\big(\Gamma_{m,l}^{\Delta}\big)^2 \big|\Delta,\Delta'\big]\,\Big] \\
\!\!\!\!\!\!\!\!\!\!\!&\leq&\!\!\!\! 
C\, \E^{\Delta,\Delta'}\!\!\left[\,\E^Y_0\left[\left(\sum_{x\in\Z^2} p_{i-1}^X(x) 
                     \E^X_0\big[L^{k-m}_{[0,i-1]}\, \big|\, F_{i-1}, X_{i-1}=x \big]\,
\E^X_{x-Y_{i-1}}\big[L_{[0,1]}(X,\Delta)\big]\right)^2 \Big| \Delta, \Delta'\right]\, \right]
\nonumber \\
\!\!\!\!\!\!\!\!\!\!\!&\leq&\!\!\!\! 
C\, \E^Y_0\E^{\Delta}\!\!\left[\!\left(\sum_{x\in\Z^2} p_{i-1}^X(x)
\Big(\E^X_0\big[L^{k-m}_{[0,i-1]}\, \big|\, F_{i-1}, X_{i-1}=x \big]\Big)^2 \! \right) \!\!\!
\left(\sum_{x\in\Z^2}p^X_{i-1}(x) \Big(\E^X_{x-Y_{i-1}}\big[L_{[0,1]}(X,\Delta)\big]\Big)^2 \!\right)\!\right], \nonumber
\end{eqnarray}
where we have applied Cauchy-Schwarz. Note that the first inner sum above does not depend on $\Delta$; while conditioned
on $Y_{i-1}$, for the second inner sum above, we have for $i\geq 2$
\begin{eqnarray}
\E^\Delta\Big[\sum_{x\in\Z^2}p^X_{i-1}(x) \Big(\E^X_{x-Y_{i-1}}\Big[L_{[0,1]}(X,\Delta)\Big]\Big)^2\Big] &\leq&
\E^\Delta\Big[\sum_{x\in\Z^2}p^X_{i-1}(x) \E^X_{x-Y_{i-1}}\Big[L_{[0,1]}(X,\Delta)\Big]\Big] \nonumber \\
&=& \int_0^1 \sum_{y\in\Z^2} p^{Y}_s(y) \sum_{x\in\Z^2} p^X_{i-1}(x) p^X_{s}(Y_{i-1}+y-x) ds \nonumber \\
&=& \int_0^1 \sum_{y\in\Z^2} p^Y_s(y) p^X_{i-1+s}(Y_{i-1}+y) ds \leq \frac{C}{i-1},
\end{eqnarray}
where we again used the local central limit theorem, and $C$ depends only on the transition kernel of $X$. Therefore, from
(\ref{A811}) we get for $i\geq 2$
\begin{eqnarray}
\E^{\Delta,\Delta'}\Big[\,\E^Y_0\big[\big(\Gamma_{m,l}^{\Delta}\big)^2 \big|\Delta,\Delta'\big]\,\Big]
&\leq& C \frac{(\log n)^{2(m-l)}}{i-1} \E^Y_0\left[\sum_{x\in\Z^2} p_{i-1}^X(x) \Big(\E^X_0\big[L^{k-m}_{[0,i-1]}\, \big|\,
F_{i-1}, X_{i-1}=x \big]\Big)^2\right] \nonumber \\
&\leq& C \frac{(\log n)^{2(k-l)}}{i-1},
\end{eqnarray}
which follows by the same calculation as in (\ref{A76}) and (\ref{A722}). Note that $\Gamma^{\Delta}_{m,l}=0$ when $i=1$.
Summing over $1\leq i\leq n$, we see that the contribution of terms involving $\Gamma^{\Delta}_{m,l}$, $l\geq 1$, to the 
variance of $\E^X_0[L_n^k|Y]$ is of order $O(\log^{2k-1} n)$. This completes the variance bound for the continuous time case.

\vspace{12pt}
\noindent{\bf Almost sure convergence of $\E^X_0[L_n^k|Y]/\log^k n$.} Because of the monotonicity of 
$\E^X_0[L^k_n |Y]$ and $\log^k n$ in $n$, we can apply the standard argument of first establishing almost sure convergence of 
$\E^X_0[L^k_n | Y]/\log^k n$ along a subsequence in $\N$ (or $\R^+$ in the continuous time case), and then use the monotonicity
to bridge the gap.

We will only treat the discrete time case. The continuous time case is identical. Fix $k\in\N$. By Theorem \ref{T:explaw}, 
$\lim_{n\to\infty} \E^{X,Y}_{0,0}[L_n^k]/\log^k n = k!/(2\pi\sqrt{\det Q})^k$, where $Q$ is its covariance matrix of
the random walk $Z=X-Y$. By the variance bound (\ref{varbound}), we have for any $\delta>0$,
\be
\P^Y_0\left(\Big| E^X_0[L_n^k|Y] -E^{X,Y}_{0,0}[L_n^k]\Big|\geq \delta (\log n)^k \right) \leq \frac {C (\log
n)^{2k-\frac{1}{2}}}{\delta^2 (\log n)^{2k}} = \frac{C}{\delta^2 (\log n)^{\frac{1}{2}}}.
\ee
Along the subsequence $t_m = e^{m^3}$, $m\in\N$, by Borel-Cantelli,
\be
\lim_{m\to\infty}\frac{\E^X_0[L_{t_m}^k|Y]}{(\log t_m)^k} = \frac{k!}{(2\pi\sqrt{\det Q})^k} \quad \text{almost surely}.
\ee
For any $t_m\leq n < t_{m+1}$, by the monotonicity of $\E^X_0[L_n^k|Y]$ and $(\log n)^k$ in $n$, we have
$$
\left(\frac{m}{m+1}\right)^{3k} \frac{\E^X_0[L_{t_m}^k|Y]}{(\log t_m)^k} 
=\frac{\E^X_0[L_{t_m}^k|Y]}{(\log t_{m+1})^k} \leq \frac{\E^X_0[L_n^k|Y]}{(\log n)^k} \leq \frac{\E^X_0[L_{t_{m+1}}^k |Y]}{(\log t_m)^k}
= \frac{\E^X_0[L_{t_{m+1}}^k |Y]}{(\log t_{m+1})^k} \left(\frac{m+1}{m}\right)^{3k}.
$$
It is then clear that $\lim_{n\to\infty} \E^X_0[L_n^k|Y]/(\log n)^k = k!/(2\pi\sqrt{\det Q})^k$ almost
surely w.r.t.\ $Y$, and Theorem \ref{T:2rw} follows by the method of moments.
\qed
\medskip

\noindent
{\bf Acknowledgement} We are grateful to Greg Lawler for showing us how Lemma \ref{L:gradpot} 
can be established for general zero mean finite variance random walks. We thank Alejandro Ram\'irez for interesting
discussions concerning the open problem formulated in the introduction, and we thank the referee for helpful 
comments. Both authors are supported by the DFG Forschergruppe 718 {\it Analysis and Stochastics 
in Complex Physical Systems}.

\end{document}